\documentclass[12pt,oneside]{amsart}
\usepackage{amsaddr}
\usepackage[utf8]{inputenc}
\usepackage[T1]{fontenc}
\usepackage{lmodern}
\usepackage{amssymb}
\usepackage{geometry}
\usepackage{amsmath}
\usepackage{amsthm}
\usepackage{graphicx}
\usepackage{color}
\usepackage[normalem]{ulem}
\usepackage[onehalfspacing]{setspace}
\usepackage{verbatim}
\usepackage[hidelinks]{hyperref}

\newcommand{\IR}{\ensuremath{\mathbb{R}}}
\newcommand{\IN}{\ensuremath{\mathbb{N}}}

\newcommand{\IP}{\ensuremath{\mathbb{P}}}
\newcommand{\IE}{\ensuremath{\mathbb{E}}}

\newcommand{\torus}{\mathbb{T}}

\newcommand{\norm}[1]{\left\Vert#1\right\Vert}
\newcommand{\set}[1]{\left\{#1\right\}}

\newcommand{\brackets}[1]{\left(#1\right)}
 
\newcommand{\scalar}[2]{\left\langle#1,#2\right\rangle}

\newcommand{\diag}{\mathop{\mathrm{diag}}}

\DeclareMathOperator{\ord}{ord}

\newtheorem{thm}{Theorem}
\newtheorem{cor}[thm]{Corollary}

\theoremstyle{plain}
\newtheorem{lemma}{Lemma}
\newtheorem{prop}{Proposition}
\theoremstyle{definition}

\newtheorem{rem}{Remark}

\title[Function values are enough for $L_2$-approximation]{Function values are enough\\ for $L_2$-approximation}

\author{
David Krieg
\and 
Mario Ullrich}
\address{Institut f\"ur Analysis, 
Johannes Kepler Universit\"at Linz, Austria}
\email{david.krieg@jku.at, mario.ullrich@jku.at}
\thanks{D.\,Krieg is supported by the Austrian Science Fund (FWF) Project F5513-N26, 
which is a part of the Special Research Program \emph{Quasi-Monte Carlo Methods:~Theory and Applications}.}
\keywords{$L_2$-approximation, sampling numbers, rate of convergence, 
random matrices, Sobolev spaces with mixed smoothness}
\subjclass[2010]{41A25, 
41A46, 
60B20; 
Secondary 
41A63, 
46E35
}

\date{\today}

\begin{document}

\begin{abstract}
We study the $L_2$-approximation of functions from a Hilbert space
and compare the sampling numbers with the approximation numbers.
The sampling number $e_n$ is the minimal worst case error that can be achieved with $n$ function values,
whereas the approximation number $a_n$ is the minimal worst case error 
that can be achieved with $n$ pieces of arbitrary linear information (like derivatives or Fourier coefficients).
We show that 
\[
 e_n \,\lesssim\, \sqrt{\frac{1}{k_n} \sum_{j\geq k_n} a_j^2},
\]
where $k_n \asymp n/\log(n)$.
This proves that the sampling numbers decay with the same polynomial rate as the approximation numbers
and therefore that function values are basically as powerful as arbitrary linear information
if the approximation numbers are square-summable.
Our result applies, in particular, 
to Sobolev spaces $H^s_{\rm mix}(\mathbb{T}^d)$ with dominating mixed smoothness 
$s>1/2$ and dimension $d\in\IN$ and we obtain
\[
e_n \,\lesssim\, n^{-s} \log^{sd}(n).
\]
For $d>2s+1$, this improves upon all previous bounds
and disproves the prevalent conjecture that Smolyak's (sparse grid) algorithm is optimal.
\end{abstract}

\maketitle

\newpage

Let $H$ be a \emph{reproducing kernel Hilbert space}, i.e., a Hilbert space of real-valued functions on a set $D$
such that point evaluation
$$
 \delta_x\colon H \to \IR,\quad \delta_x(f) := f(x)
$$
is a continuous functional for all $x\in D$.
We consider numerical approximation of functions from such spaces, 
using only function values.
We measure the error in the space $L_2=L_2(D,\mathcal{A},\mu)$
of square-integrable functions with respect to an arbitrary 
measure $\mu$ such that $H$ is embedded into $L_2$.
This means that the functions in $H$ are square-integrable
and two functions from $H$ that are equal \mbox{$\mu$-almost} everywhere
are also equal point-wise.

We are interested in the \emph{$n$-th minimal worst-case error}
\[
e_n \,:=\, e_n(H) \,:=\, 
\inf_{\substack{x_1,\dots,x_n\in D\\ \varphi_1,\dots,\varphi_n\in L_2}}\, 
\sup_{f\in H\colon \|f\|_H\le1}\, 
\Big\|f - \sum_{i=1}^n f(x_i)\, \varphi_i\Big\|_{L_2},
\]
which is the worst-case error of an optimal algorithm that  
uses at most $n$ function values.
These numbers are sometimes called \emph{sampling numbers}. 
We want to compare $e_n$ with the \emph{$n$-th approximation number} 
\[
a_n \,:=\, a_n(H) \,:=\,
\inf_{\substack{L_1,\dots,L_n\in H'\\ \varphi_1,\dots,\varphi_n\in L_2}}\, 
\sup_{f\in H\colon \|f\|_H\le1}\, 
\Big\|f - \sum_{i=1}^n L_i(f)\, \varphi_i\Big\|_{L_2},
\]
where $H'$ is the space of all bounded, linear functionals on $H$.
This is the worst-case error of an optimal algorithm that uses 
at most $n$ linear functionals as information.
Clearly, we have $a_n\leq e_n$ since the point evaluations form a subset of $H'$.

The approximation numbers are quite well understood in many cases because 
they are equal to the singular values of the embedding operator ${\rm id}\colon H\to L_2$.
However, the sampling numbers still resist a precise analysis. 
For an exposition of such approximation problems we refer to 
\cite{NW08,NW10,NW12}, especially~\cite[Chapter~26~\&~29]{NW12}, 
and references therein.
One of the fundamental questions in the area asks for the relation 
of $e_n$ and $a_n$ for specific Hilbert spaces $H$. 
The minimal assumption on $H$ is the compactness of the embedding ${\rm id}\colon H\to L_2$.
It is known that
\[
 \lim_{n\to\infty} e_n = 0
 \quad \Leftrightarrow \quad
 \lim_{n\to\infty} a_n = 0 \quad
 \quad \Leftrightarrow \quad
 H \hookrightarrow L_2 \text{ compactly},
\]
see \cite[Section~26.2]{NW12}. 
However, the compactness of the embedding is not enough for
a reasonable comparison of the speed of this convergence, see \cite{HNV08}.
If $(a_n^*)$ and $(e_n^*)$ are decreasing sequences that converge to zero
and $(a_n^*)\not\in\ell_2$, one may construct $H$ and $L_2$ such that
$a_n=a_n^*$ for all $n\in\IN$
and $e_n\geq e_n^*$ for infinitely many $n\in\IN$.
In particular, if
\[
 \ord(c_n) := \sup\set{s\geq 0 \colon \lim_{n\to\infty} c_n n^{s}=0}
\]
denotes the (polynomial) order of convergence of a positive sequence $(c_n)$,
it may happen that $\ord(e_n)=0$ even if $\ord(a_n)=1/2$.

It thus seems necessary to assume that $(a_n)$ is in $\ell_2$, i.e.,
that ${\rm id}\colon H\to L_2$ is a Hilbert-Schmidt operator.
This is fulfilled, e.g., for \emph{Sobolev spaces} defined on the unit cube, see Corollary~\ref{cor:sob}.
Under this assumption, it is proven in \cite{KWW09} that
\[
 \ord(e_n) \geq \frac{2 \ord(a_n)}{2\ord(a_n) + 1}\, \ord(a_n).
\]
In fact, the authors of \cite{KWW09} conjecture that the order of convergence is the same for both sequences. 
We give an affirmative answer to this question.
Our main result can be stated as follows.

\begin{thm}\label{thm:main}
 There are absolute constants $C,c>0$ and a sequence of natural numbers
 $(k_n)$ with $k_n\ge c n/\log(n+1)$ such that the following holds.
 For any~$n\in\IN$, any measure space $(D,\mathcal A,\mu)$ and 
 any reproducing kernel Hilbert space $H$ of real-valued functions on $D$
 that is embedded into $L_2(D,\mathcal A,\mu)$, we have
 \[
  e_n(H)^2 \,\le\, \frac{C}{k_n} \sum_{j\geq k_n} a_j(H)^2.
 \]
\end{thm}

In particular, we obtain the following result on the order of convergence.
This solves Open Problem~126 in~\cite[p.~333]{NW12},
see also~\cite[Open Problems 140~\&~141]{NW12}.

\begin{cor}\label{cor:order-of-convergence}
 Consider the setting of Theorem~\ref{thm:main}.
 If $a_n(H)\lesssim n^{-s}\log^\alpha(n)$ for some $s>1/2$ and $\alpha\in\IR$, 
 then we obtain
\[
 e_n(H) \,\lesssim\, n^{-s}\log^{\alpha+s}(n).
\]
 In particular, we always have $\ord(e_n)=\ord(a_n)$.
\end{cor}

Let us now consider a specific example. 
Namely, we consider \emph{Sobolev spaces with (dominating) mixed smoothness} 
defined on the $d$-dimensional torus $\mathbb{T}^d \cong [0,1)^d$. 
These spaces attracted quite a lot of attention in various areas of 
mathematics due to their intriguing attributes in high-dimensions. 
For history and the state of the art 
(from a numerical analysis point of view) 
see~\cite{DTU16,Tem18,Tri10}.

Let us first define a one-dimensional and real-valued orthonormal basis 
of $L_2(\torus)$ by
$b_0^{(1)}=1$, $b_{2m}^{(1)}=\sqrt{2}\cos(2\pi m x)$ and 
$b_{2m-1}^{(1)}=\sqrt{2}\sin(2\pi m x)$
for $m\in\IN$.
From this we define a basis of $L_2(\torus^d)$ using 
$d$-fold tensor products: We set
$\mathbf{b}_{\bf m}:=\bigotimes_{j=1}^d b_{m_j}^{(1)}$ 
for ${\bf m}=(m_1,\dots,m_d)\in\IN_0^d$.
The Sobolev space with dominating mixed smoothness $s>0$ can be defined as
\[
 H :=
 H^s_{\rm mix}(\torus^d)
 := \Big\{ f \in L_2(\torus^d)
 \,\Big|\, \|f\|_H^2 :=  \sum_{{\bf m} \in \IN_0^d} \prod_{j=1}^d(1+|{m_j}|^{2s}) \scalar{f}{\bf b_{{\bf m}}}_{L_2}^2 <\infty \Big\}.
\]
This is a Hilbert space. It satisfies our assumptions
whenever $s>1/2$. 
It is not hard to prove that an equivalent norm in $H^s_{\rm mix}(\torus^d)$ 
for $s\in\IN$ is given by
\[
\|f\|_{H^s_{\rm mix}(\torus^d)}^2 
\,=\, \sum_{\alpha\in\{0,s\}^d} \|D^\alpha f\|_{L_2}^2.
\]
The approximation numbers $a_n = a_n(H)$ 
are known for some time to satisfy
\[
a_n \,\asymp\, n^{-s} \log^{s(d-1)}(n)
\]
for all $s>0$, see e.g.~\cite[Theorem~4.13]{DTU16}.
The sampling numbers $e_n = e_n(H)$, 
however, seem to be harder to tackle. The best bounds so far are 
\[
n^{-s} \log^{s(d-1)}(n) \,\lesssim\; e_n 
\;\lesssim\, n^{-s} \log^{(s+1/2)(d-1)}(n)
\]
for $s>1/2$. 
The lower bound easily follows from $e_n\ge a_n$, 
and the upper bound was proven in \cite{SU07},
see also \cite[Chapter~5]{DTU16}.
For earlier results on this prominent problem, see 
\cite{Si03,Si06,Tem93,Ul08}.
Note that finding the right order of $e_n$ in this case is posed as 
\emph{Outstanding Open Problem~1.4} in~\cite{DTU16}.
From Corollary~\ref{cor:order-of-convergence}, setting $\alpha=s(d-1)$ in the second part, 
we easily obtain the following.

\begin{cor}\label{cor:sob}
Let $H^s_{\rm mix}(\torus^d)$ be the Sobolev space with mixed smoothness
as defined above. Then, for $s>1/2$, we have
\[
e_n\big(H^s_{\rm mix}(\torus^d)\big) \,\lesssim\,  n^{-s} \log^{sd}(n).
\]
\end{cor}

The bound in Corollary~\ref{cor:sob} improves on the previous bounds 
if $d>2s+1$, or equivalently $s<(d-1)/2$. 
With this, we disprove Conjecture~5.26 from~\cite{DTU16} and show, 
in particular, that Smolyak's algorithm is not optimal in these cases.
Although our techniques do not lead to 
an explicit deterministic 
algorithm that achieves the above bounds, 
it is interesting that $n$ i.i.d.~random points 
are suitable with positive probability (independent of $n$).
\medskip

Let us conclude with a few topics for future research.
While this paper was under review, 
Theorem~\ref{thm:main} has already been extended
to the case of complex-valued functions and 
non-injective operators ${\rm id}\colon H\to L_2$ in~\cite{KUV19},
including explicit values for the constants $c$ and $C$, 
see also~\cite{U20}. 
It remains open to generalize our results to non-Hilbert space settings.
It is also quite a different question whether the sampling numbers 
and the approximation numbers 
behave similarly with respect to the dimension of the domain~$D$.
This is a subject of tractability studies.
We refer to \cite[Chapter~26]{NW12} and especially \cite[Corollary~8]{NW16}. 
Here, we only note that the constants of Theorem~\ref{thm:main} are, 
in particular, independent of the domain, and that this may
be utilized for these studies, see also~\cite{KUV19}.

\goodbreak

\section*{The Proof}

The result follows from a combination of the general technique to 
assess the quality of \emph{random information}
as developed in~\cite{HKNPU19b,HKNPU19a},  
together with bounds on the singular values of random matrices with 
independent rows from~\cite{MP06}. 

Before we consider algorithms that only use function values, let us briefly recall the situation for arbitrary linear functionals. 
In this case, the minimal worst-case error $a_n$ is given via the singular value decomposition of ${\rm id}: H\to L_2$ in the following way.
Since $W={\rm id}^*{\rm id}$ is positive, compact and injective, 
there is an orthogonal basis $\mathcal B=\set{b_j \colon j\in\IN}$ of $H$ that consists of eigenfunctions of $W$.
Without loss of generality, we may assume that $H$ is infinite-dimensional. 
It is easy to verify that $\mathcal B$ is also orthogonal in $L_2$.
We may assume that the eigenfunctions are normalized in $L_2$
and that $\|b_1\|_H \leq \|b_2\|_H \leq \dots$.
From these properties, it is clear that the Fourier series
\[
f\,=\,\sum_{j=1}^\infty f_j b_j, 
\qquad \text{ where } \quad f_j:=\scalar{f}{b_j}_{L_2},
\]
converges in $H$ for every $f\in H$,
and therefore also point-wise.
The optimal algorithm based on $n$ linear functionals is given by
\[
 P_n: H \to L_2, \quad P_n(f):=\sum_{j\leq n} f_j b_j,
\]
which is the $L_2$-orthogonal projection onto $V_n:={\rm span}\{b_1,\hdots,b_n\}$.
We refer to \cite[Section~4.2]{NW08} for details. 
We obtain that
\[
 a_n(H)=\sup_{f\in H\colon \|f\|_H\le1} \big\Vert f - P_n(f) \big\Vert_{L_2} =\|b_{n+1}\|_H^{-1}.
\]

We now turn to algorithms using only function values.
In order to bound the minimal worst-case error $e_n$ from above, we 
employ the \emph{probabilistic method} in the following way.
Let $x_1,\dots,x_n\in D$ be i.i.d.~random variables with $\mu$-density
\[
 \varrho: D\to \IR, \quad \varrho(x) := \frac12 \left(
 \frac1k \sum_{j< k} b_{j+1}(x)^2  +  \frac{1}{\sum_{j\geq k} a_j^2} \sum_{j\geq k} a_j^2 b_{j+1}(x)^2
 \right),
\]
where $k\leq n$ will be specified later.
Given these sampling points, we consider the algorithm 
\[
 A_n: H\to L_2, \quad A_n(f):=\sum_{j=1}^k (G^+ N f)_j b_j,
\]
where $N\colon H\to \IR^n$ with $N(f):=(\varrho(x_i)^{-1/2}f(x_i))_{i\leq n}$ 
is the weighted \emph{information mapping} and
$G^+\in \IR^{k\times n}$ is the Moore-Penrose inverse of the matrix 
\[
 G:=(\varrho(x_i)^{-1/2} b_j(x_i))_{i\leq n, j\leq k} \in \IR^{n\times k}.
\]
This algorithm is a weighted least squares estimator:
If $G$ has full rank, then 
\[
 A_n(f)=\underset{g\in V_k}{\rm argmin}\, \sum_{i=1}^n \frac{\vert g(x_i) - f(x_i) \vert^2}{\varrho(x_i)}. 
\]
In particular, we have $A_n(f)=f$ whenever $f\in V_k$.
The \emph{worst-case error} of $A_n$ is defined as
\[
 e(A_n) := \sup_{f\in H\colon \|f\|_H\le1}\,  \big\|f - A_n(f)\big\|_{L_2}.
\]
Clearly, we have $e_n\leq e(A_n)$ for every realization of $x_1,\dots,x_n$.
Thus, it is enough to show that $e(A_n)$ obeys the desired upper bound with positive probability. 

\begin{rem}
If $\mu$ is a probability measure and if the basis is uniformly bounded, i.e.,
if $\sup_{j\in\IN}\, \Vert b_j\Vert_\infty < \infty$,
we may also choose $\varrho\equiv 1$ and consider i.i.d.~sampling points with distribution~$\mu$.
\end{rem}

\begin{rem}
Weighted least squares estimators are widely studied in the literature. 
We refer to \cite{Bj96,CM17}.
In contrast to previous work, we show that we can choose a fixed set of weights and sampling points that work \emph{simultaneously} for all $f\in H$.
We do not need additional assumptions on the function $f$, the basis $(b_j)$ or the measure $\mu$.
For this, we think that our modification of the weights is important.
\end{rem}

\begin{rem}
The worst-case error $e(A_n)$ of the randomly chosen algorithm $A_n$
is not to be confused with the Monte Carlo error of a randomized algorithm,
which can be defined by
\[
 e^{\rm ran}(A_n) \,:=\, \sup_{f\in H\colon \|f\|_H\le1}\, 
\left(\IE \left\|f - A_n(f)\right\|_{L_2}^2 \right)^{1/2}.
\]
The Monte-Carlo error is a weaker error criterion.
It is shown in \cite{Kr19}, see also \cite{WW06}, 
that the assumptions of Corollary~\ref{cor:order-of-convergence}
give rise to a randomized algorithm $M_n$
which uses at most $n$ function values and satisfies
\[
 e^{\rm ran}(M_n) \,\lesssim\, n^{-s}\log^\alpha(n).
\]
However, this does not imply that the worst-case error 
$e(M_n)$ is small for any realization of~$M_n$.
\end{rem}

To give an upper bound on $e(A_n)$,
let us assume that $G$ has full rank. 
For any $f\in H$ with $\Vert f\Vert_H\leq 1$, we have
\[\begin{split}
 \norm{f-A_n f}_{L_2} \,&\le\, a_k + \norm{P_k  f - A_n f}_{L_2}
 \,=\, a_k + \norm{A_n(f- P_k f)}_{L_2} \\
 &=\, a_k + \norm{G^+ N(f- P_k f)}_{\ell_2^k} \\
 &\le\, a_k +\norm{G^+\colon \ell_2^n \to \ell_2^k} \norm{N\colon P_k(H)^\perp \to \ell_2^n}.
\end{split}
\]
The norm of $G^+$ is the inverse of the $k$th largest (and therefore the smallest) singular value of the matrix $G$.
The norm of $N$ is the largest singular value of the matrix 
\[
 \Gamma :=\big(\varrho(x_i)^{-1/2} a_j b_{j+1}(x_i)  \big)_{1\leq i \leq n, j\geq k} \in \IR^{n\times \infty}.
\]
To see this, note that $N=\Gamma \Delta$ on $P_k(H)^\perp$, 
where the mapping $\Delta\colon P_k(H)^\perp\to \ell_2$ 
with $\Delta g=(g_{j+1}/a_j)_{j\ge k}$
is an isomorphism.
This yields 
\begin{equation}
\label{eq:basic}
 e(A_n) \leq a_k + \frac{s_{\rm max}(\Gamma)}{s_{\rm min}(G)}.
\end{equation}
It remains to bound $s_{\rm min}(G)$ from below and $s_{\rm max}(\Gamma)$ from above.
Clearly, any nontrivial lower bound on $s_{\rm min}(G)$ automatically yields that the matrix $G$ has full rank.
To state our results, let
\[
\beta_k 
\,:=\, \brackets{\frac{1}{k} \sum_{j\geq k} a_j^2}^{1/2} 
\qquad\text{ and }\qquad
\gamma_{k}\,:=\,\max\Big\{a_k,\,\beta_k\Big\}.
\]
Note that $a_{2k}^2\le\frac1k(a_k^2+\hdots+a_{2k}^2)\le \beta_{k}^2$\, 
for all $k\in\IN$ 
and thus $\gamma_{k} \leq \beta_ {\lfloor k/2 \rfloor}$.
Before we continue with the proof of Theorem~\ref{thm:main}, 
we show that Corollary~\ref{cor:order-of-convergence} follows from 
Theorem~\ref{thm:main} by providing the order of $\beta_k$ in the following special case. 
The proof is an easy exercise.
\begin{lemma}\label{lem:beta}
Let $a_n\asymp n^{-s}\log^{\alpha}(n)$ for some $s,\alpha\in\IR$. Then, 
\[
\beta_{k} \,\asymp\, \begin{cases}
a_k, & \text{if } s>1/2, \\
a_k \sqrt{\log(k)}, & \text{if } s=1/2 \,\text{ and }\, \alpha<-1/2,\\
\end{cases}
\]
and $\beta_k=\infty$ in all other cases.
\end{lemma}


The rest of the paper is devoted to the proof of the following two claims: 
There exist constants $c,C>0$ such that, 
for all $n\in \IN$ and $k= \lfloor c\,n/\log n\rfloor$,
we have \smallskip\\[1mm]
\noindent{\bf Claim 1:} 
\[
\IP\Big(s_{\rm max}(\Gamma) \,\leq\, C\, \gamma_{k}\, n^{1/2} \Big) > 1/2.
\]

\smallskip
\noindent{\bf Claim 2:} 
\[
\IP\Big(s_{\rm min}(G) \,\geq\, n^{1/2}/2 \Big) > 1/2.
\medskip
\]

Together with \eqref{eq:basic}, this will yield with positive probability that
\[
 e(A_n) \,\le\, a_k + 2C\,\gamma_{k} 
 \leq (2C+1)\, \gamma_{k} \leq (2C+1)\, \beta_{\lfloor k/2 \rfloor},
\]
which is the statement of Theorem~\ref{thm:main}.

\medskip
Both claims are based on~\cite[Theorem~2.1]{MP06}, which we state here 
in a special case.
Recall that, for $X\in \ell_2$, the operator $X\otimes X$
is defined on $\ell_2$ by $X\otimes X(v)= \langle X,v\rangle_2\cdot X$.
By $\norm{M}$ we denote the spectral norm of a matrix $M$.

\begin{prop}\label{prop:MP}
There exists an absolute constant $c>0$ for which the following holds. 
Let $X$ be a random vector in $\IR^k$ or $\ell_2$ with  
$\|X\|_2\le R$ with probability 1, and let $X_1,X_2,\dots$ 
be independent copies of $X$. 
We put
\[
D:=\IE(X\otimes X),
\qquad
A \,:=\, R^2\, \frac{\log n}{n}
\qquad\text{ and }\qquad
B \,:=\, R\, \|D\|^{1/2} \sqrt{\frac{\log n}{n}}. 
\]
Then, for any $t>0$, 
\[
\IP\left(\bigg\|\sum_{i=1}^n X_i\otimes X_i - nD\bigg\|
	\,\ge\, c\, t\, \max\{A, B\}\, n\right)
\,\le\, 2e^{-t}.
\]
\end{prop}

\begin{proof}[{\bf Proof of Proposition~1}]
We describe the steps needed to obtain this reformulation 
of~\cite[Theorem~2.1]{MP06}.
For this let  
$\|Z\|_{\psi_\alpha}:=\inf\{C>0\colon\IE\exp(|Z|^{\alpha}/C^{\alpha})\le2\}$
for $Z=\|X\|_2$
and 
$$
\rho:=\sup\left\{\left(\IE\scalar{X}{\theta}_2^4\right)^{1/4}\colon 
	{\theta\in\IR^k \text{ with } \|\theta\|_2=1}\right\}.
$$
Theorem~2.1 of~\cite{MP06} then states that  
\[
\IP\left(\bigg\|\sum_{i=1}^n X_i\otimes X_i - nD\bigg\|
	\,\ge\, c\, t\, \max\{\widetilde{A}, \widetilde{B}\}\, n\right)
\,\le\, 
2e^{-t^{\alpha/(2+\alpha)}}
\]
with
\[
\widetilde{A} \,:=\, \|Z\|_{\psi_\alpha}^2\frac{(\log n)^{1+\frac{2}{\alpha}}
}{n}
\qquad\text{ and }\qquad
\widetilde{B} \,:=\, \frac{\rho^2}{\sqrt{n}} + \left(\|D\|\cdot \widetilde{A} \right)^{1/2}, 
\]
for all $t>0$, $\alpha\ge1$, $n\in\IN$ and some absolute constant $c>0$.
Note that the 2 in the right hand side of above inequality is 
missing in \cite[Theorem~2.1]{MP06}, but can be found in the proof.

From $\|X\|_2\le R$ we obtain $\|Z\|_{\psi_\alpha}\le2R$ for all $\alpha\ge1$. 
Therefore, we can take the limit $\alpha\to\infty$ and obtain the 
result with $\widetilde{A}=\frac{R^2 \log n}{n}$ and 
$\widetilde{B}=\frac{\rho^2}{\sqrt{n}} + R\|D\|^{1/2}\sqrt{\frac{\log n}{n}}$ 
(and a slightly changed constant $c$). 
Moreover, we have 
\[
\IE\scalar{X}{\theta}_2^4
\le R^2 
 \cdot \IE\scalar{X}{\theta}_2^2
= R^2 \cdot \scalar{D\theta}{\theta}_2 
\le R^2 \cdot \|D\|
\] 
for any $\theta\in\IR^k$ (or $\ell_2$) with $\|\theta\|_2=1$, 
which implies $\rho^2\le R\cdot\|D\|^{1/2}$.
This ``trick'' leads to an improvement over~\cite[Corollary~2.6]{MP06}
and yields 
our formulation of the result. \\
\end{proof}

\medskip

\begin{proof}[{\bf Proof of Claim 1}]
Consider independent copies $X_1,\hdots,X_n$ of the vector 
\[
X:=\varrho(x)^{-1/2} (a_k b_{k+1}(x), a_{k+1} b_{k+2}(x), \hdots),
\]
where $x$ is a random variable on $D$ with density $\varrho$.
Clearly, $\sum_{i=1}^n X_i\otimes X_i = \Gamma^* \Gamma$ with $\Gamma$ 
from above.
First observe
\[
 \norm{X}_2^2 \,=\, \varrho(x)^{-1} \sum_{j\geq k} a_j^2\, b_{j+1}(x)^2 
 \,\leq\, 2 \sum_{j\geq k} a_j^2
 \,=\, 2 k\, \beta_{k}^2 \,=:\, R^2.
\]
Since $D:=\IE(X\otimes X)=\diag(a_k^2, a_{k+1}^2, \hdots)$ 
we have $\|D\|=a_k^2$. 
This implies, with $A$ and $B$ defined as in 
Proposition~\ref{prop:MP}, that
\[
A \,\le\, 2 k\, \beta_{k}^2\, \frac{\log n}{n}
\]
and 
\[
B \,\le\, (2 k\, \beta_{k}^2\,)^{1/2} a_k\, \sqrt{\frac{\log n}{n}}.
\]
Choosing $k= \lfloor c\,n/\log n\rfloor$ for $c$ small enough, 
we obtain
\[
 \IP\Big(\norm{\Gamma^*\Gamma - nD} \geq t\,\gamma_{k}^2\, n\Big)
 \leq 2\exp\brackets{-t}.
\]
By choosing $t=2$, we obtain with probability greater $1/2$ that
\[
 s_{\rm max}(\Gamma)^2 = \norm{\Gamma^*\Gamma} \leq \norm{nD} + \norm{\Gamma^*\Gamma - nD}
 \leq n\, a_k^2 + 2 \gamma_{k}^2 n
 \leq 3\, \gamma_{k}^2\, n.
\]
This yields Claim~1.\\
\end{proof}

\begin{proof}[{\bf Proof of Claim 2}]
Consider $X:=\varrho(x)^{-1/2}(b_1(x), \hdots, b_k(x))$
with $x$ distributed according to $\varrho$.
Clearly, $\sum_{i=1}^n X_i\otimes X_i = G^*G$ with $G$ from above.
First observe
\[
 \norm{X}_2^2 \,=\, \varrho(x)^{-1} \sum_{j\le k} b_j(x)^2 
 \,\leq\, 2 k \,=:\, R^2.
\]
Since $D:=\IE(X\otimes X)=\diag(1, \hdots,1)$ 
we have $\|D\|=1$. 
This implies, with $A$ and $B$ defined as in 
Proposition~\ref{prop:MP}, that
\[
A \,\le\, 2 k\, \frac{\log n}{n}
\]
and 
\[
B \,\le\, (2 k)^{1/2} \sqrt{\frac{\log n}{n}}.
\]
Again, choosing $k= \lfloor c\,n/\log n\rfloor$ for $c$ small enough, 
we obtain
\[
 \IP\left(\norm{G^*G - nD} \geq \frac{t\, n}{4}\right)
 \leq 2\exp\brackets{-t}.
\]
By choosing $t=2$, we obtain with probability greater $1/2$ that
\[
 s_{\rm min}(G)^2 = s_{\rm min}(G^*G) \,\geq\, s_{\rm min}(nD) - \|G^*G - nD\|
 \,\geq\, 
 n/2.
\]
This yields Claim~2.\\
\end{proof}


\end{document}